\documentclass[12pt, a4paper]{article}
\usepackage[T1]{fontenc}
    \usepackage{lmodern}
\usepackage[utf8]{inputenc}

\usepackage{subcaption}
\usepackage{amssymb}
\usepackage{mathtools}
\usepackage{amsthm}
\usepackage{microtype}
\usepackage{hyphenat}
\usepackage{enumitem}
\usepackage{xcolor}
\usepackage{mdframed}
\usepackage{pdfpages}
\usepackage{csquotes}
\usepackage{tikz-cd}
	\usetikzlibrary{arrows}
	\tikzcdset{arrow style=tikz, diagrams={>=angle 90}}
\usepackage{booktabs}
    \usetikzlibrary{calc}
\usepackage[hidelinks]{hyperref}
\hypersetup{colorlinks = true, linkcolor = blue, citecolor= blue, urlcolor= gray}
\usepackage{standalone}

\DeclareMathOperator{\ban}{ban}
\DeclareMathOperator{\Int}{int}

\DeclarePairedDelimiter{\abs}{\lvert}{\rvert}

\DeclarePairedDelimiter{\set}{\{}{\}}
\DeclarePairedDelimiter{\seq}{(}{)}
\DeclarePairedDelimiterX{\oo}[2]{(}{)}{#1,#2}
\DeclarePairedDelimiterX{\oc}[2]{(}{]}{#1,#2}
\DeclarePairedDelimiterX{\co}[2]{[}{)}{#1,#2}
\DeclarePairedDelimiterX{\cc}[2]{[}{]}{#1,#2}
\DeclarePairedDelimiterX{\ip}[2]{\langle}{\rangle}{#1 , #2}
\DeclarePairedDelimiterX{\euler}[2]{\langle}{\rangle}{#1,#2}
\DeclarePairedDelimiterX{\symbi}[2]{(}{)}{#1,#2}
\DeclarePairedDelimiterX{\pres}[2]{\langle}{\rangle}{#1 \:\vert\: #2}

\newcommand{\s}{\sigma}
\newcommand{\NN}{\mathbb{N}}

\newcommand{\RR}{\mathbb{R}}

\newcommand{\mcL}{\mathcal{L}}
\newcommand{\mcA}{\mathcal{A}}

\newcommand{\mcR}{\mathcal{R}}

\newcommand{\mcW}{\mathcal{W}}

\newcommand{\tops}{\mathrm{top}}
\newcommand{\G}{\mathcal{G}}
\newcommand{\F}{\mathcal{F}}
\newcommand{\V}{\mathcal{V}}

\theoremstyle{plain}
\newtheorem{thm}{Theorem}[section]
\newtheorem{lem}[thm]{Lemma}
\newtheorem{prop}[thm]{Proposition}
\newtheorem{cor}[thm]{Corollary}

\newtheorem*{conjecture}{Conjecture}

\theoremstyle{definition}
\newtheorem{defn}[thm]{Definition}
\newtheorem{ex}[thm]{Example}

\theoremstyle{remark}

\newtheorem*{rmk}{Remark}

\newtheorem*{note}{Note}

\allowdisplaybreaks[4]

\title{
Graph Iterated Function Systems
and Fractal Tops}
\author{Grover Lancaster-Cole \\ Georgiana Lyall \\ Thomas Malcolm \\ Qiyu Zhou}
\date{}

\begin{document}

\maketitle

\begin{abstract}
Following the work of Louisa and Michael Barnsley \cite{blowups_tops} on results in tops of iterated function systems, we extend their work to graph-directed iterated function systems by investigating the relationship between top addresses and shift spaces. For the simplest overlapping interval IFS, we find a sufficient condition for its tops code space's closure to be a shift space of finite type. Likewise, we find that shift invariance properties do not directly extend to the graph-directed setting.
\end{abstract}

\section{Introduction}


The notion of a \emph{graph iterated function system} (GIFS) was introduced in \cite{gifs} as a generalisation of the iterated function system (IFS), a system of maps used to generate fractal images \cite{fractals_everywhere}. The notion of a \emph{fractal top} for a classical iterated function system was introduced in \cite{blowups_tops} with the primary purpose of constructing fractal tilings of the plane. In this paper, we study fractal tops for graph iterated function systems, with a particular focus on their relationships to \emph{shift spaces} and shift invariance properties.

\section{Shift spaces}

The material in this section is mostly adapted from \cite{coding}. However, where they study bi-infinite strings, we shall study right-infinite strings, that is, functions from the natural numbers $\NN$. We follow the convention that $0\notin \NN$, and we write $[N]=\set{1,\dots,N}$.

Let $\mcA$ be a finite set of symbols. A \emph{word} $\alpha$ is a finite string of elements of $\mcA$, we denote its length by $|\alpha|$ and the collection of all words in $\mcA$ by $\mcW$.
For a subset $B$ of $\mcW$, we denote by $\ban(B)$ the collection of strings in $\mcA^\NN$ which do not contain any of the words in $B$.

\begin{defn}
A \emph{shift space} over $\mcA$ is a subset $X$ of $\mcA^\NN$ such that $X=\ban(B)$ for some collection $B$ of words.
A shift space is said to be of \emph{finite type} if $B$ can be chosen to be finite.
\end{defn}

Given $w,w'\in\mcW$, we may form the concatenation $ww'\in\mcW$ and, similarly, given $w\in\mcW$ and $x\in \mcA^\NN$, we may form the concatenation $wx\in \mcA^\NN$. We have a metric $d$ on $\mcA^\NN$ defined by
\[ d(x,y) = 2^{-\min\set{k:x_k\ne y_k}} \]
and a continuous \emph{shift operator} $S:\mcA^\NN\to \mcA^\NN$ given by
\[S(x_1x_2x_3\cdots)=x_2x_3x_4\cdots. \]
For $x\in \mcA^\NN$ where $x=x_1x_2x_3\cdots $, we write $x|_{[i, j]} = x_i\cdots x_j$ for $1\leq i <j$.
It is clear that a sequence $\seq{x_n}_{n\in\NN}$ converges to $x$ in $\mcA^\NN$ if and only if for every $k\in\NN$ there exists $N\in\NN$ such that $x_n|_{[1,k]} = x|_{[1,k]}$ whenever $n\ge N$.
It is straightforward to verify that $\mcA^\NN$ is compact.

For any subset $A$ of $\mcA^\NN$, the \emph{language} $\mcL_A$ associated with $A$ consists of the words $w$ which appear in some $x\in A$.
Words in $\mcL_A^c$ are said to be \emph{banned}, and a banned word $w$ is said to be \emph{reduced} if $w$ has no proper banned subwords. The set of reduced banned words will be denoted by $\mcR_A$. It is easy to check that every banned word has a reduced banned subword, and that $\ban(\mcR_A)=\ban(\mcL_A^c)$.

\begin{prop}\label{prop:character}
For a subset $X$ of $\mcA^\NN$, the following are equivalent:
\begin{enumerate}
    \item $X$ is a shift space
    \item $X=\ban(\mcR_X)$
    \item $X$ is closed and $S(X)\subseteq X$.
\end{enumerate}
\end{prop}

\begin{proof}
Let $B\subseteq \mcW$ and let $(x_n)_{n\in\NN}$ be a sequence in $\ban(B)$ converging to some $x\in \mcA^\NN$. If $w$ is a word appearing in $x$, then $w$ must appear in some $x_n$, so $w\notin B$. So $x\in \ban(B)$, showing that $\ban(B)$ is closed.
Now, let $x\in \ban(B)$. If a word $w$ appears in $S(x)$ then $w$ appears in $x$, so $w\notin B$. That is, $S(x)\in \ban(B)$, so $S(\ban(B))\subseteq \ban(B)$. We have shown that (a) implies (c).

Now, let $X\subseteq\mcA^\NN$ be closed and suppose $S(X)\subseteq X$. It is clear that $X\subseteq \ban(\mcL_X^c)$.
For the opposite inclusion, let $x\in \ban(\mcL_X^c)$. Given $n\in\NN$, we have $x|_{[1,n]}\in \mcL_X$. But, since $S^k(X)\subseteq X$ for every $k$, in fact $x|_{[1,n]}y\in X$ for some $y\in X$. Since $n$ was arbitrary and $X$ is closed, we conclude that $x\in X$. Hence $\ban(\mcL_X^c)\subseteq X$, that is, $X=\ban(\mcL_X^c)=\ban(\mcR_X)$. So (c) implies (b).

Clearly (b) implies (a).
\end{proof}

\begin{rmk}
It is not true in general that $S(X)=X$ when $X$ is a shift space.
For a simple counterexample, take $B=\set{12,22}$. Then $2\overline 1 \in \ban(B)$, but $S^{-1}(2\overline 1)$ consists of $12\overline 1$ and $22\overline 1$, neither of which are in $\ban(B)$.
This is a key difference between our setting and that in \cite{coding}. A subset $A$ of $\mcA^\NN$ is said to be \emph{shift invariant} if $S(A)=A$.
\end{rmk}

From here on we shall assume that $\mcA$ is equipped with a total order, and we endow $\mcA^\NN$ with the associated lexicographic order.

\begin{lem}\label{lem:closed}
If $K$ is a nonempty closed subset of $\mcA^\NN$, then $K$ has a minimal element.
\end{lem}

\begin{proof}
Let $K$ be a nonempty closed subset of $\mcA^\NN$.
Pick $x_1\in K$ such that $(x_1)_1$ is the minimal $i\in \mcA$ such that there exists $y\in \mcA^\NN$ with $iy \in K$. For each $n\in\NN$, pick $x_{n+1}\in K$ such that $x_{n+1}|_{[1,n]} = x_n|_{[1,n]}$ and $(x_{n+1})_{n+1}$ is the minimal $i\in \mcA$ such that there exists $y\in \mcA^\NN$ with $x_n|_{[1,n]}i y \in K$. This defines a Cauchy sequence $\seq{x_n}_{n\in\NN}$, which converges to some $x$ in $K$ because $K$ is complete. It is clear from this construction that $x$ must be a minimal element of $K$.
\end{proof}

We will denote the minimal element of $K$ by $\min K$. The next result is a translation of a Lemma proved in \cite{blowups_tops}.

\begin{lem}\label{lem:minimum}
Let $K$ be a nonempty closed subset of a shift space $X$ and let $x=\min K$. Then $\set{y \in X : x_1 y\in K}$ is nonempty and closed, and
\[ x = x_1 \min\set{y \in X : x_1 y\in K}.\]
\end{lem}

\begin{proof}
Let $L=\set{y\in X:x_1y\in K}$.
Suppose $\seq{y_n}_{n\in\NN}$ is a sequence in $L$ which converges to some $y$ in $X$. Then $\seq{x_1y_n}_{n\in\NN}$ converges to $x_1y$ in $K$, so $y\in L$. This shows that $L$ is closed.
Clearly $L$ is nonempty because it contains $S(x)$. In particular, $x_1\min L \le x_1 S(x) = x$.
But $x_1\min L\in K$, so that $x\le x_1\min L$. Hence $x = x_1 \min L$.
\end{proof}

We end this section with the note that $\mcW$ has the structure of a partial order given by writing $w<w'$ if $wx\le w'x'$ whenever $x,x'\in \mcA^\NN$.

\section{Graph-directed iterated function systems}

This brief section is devoted to providing the basic definitions and notation for graph-directed iterated function systems, which are the primary objects of interest for this paper. The reader is directed to \cite{tiling_ifs} for a comprehensive introduction to this theory.






\begin{defn}
A \emph{graph iterated function system} (graph IFS) on $\RR^M$ is a pair $(\F,\G)$, where $\F=\set{f_1,\dots,f_N}$ is a finite collection of invertible contractive maps on $\RR^M$ and $\G$ is a strongly connected primitive directed graph with $N$ edges labelled $1,\dots,N$.
\end{defn}




We denote by $i^-$ and $i^+$ the source and target vertex for the edge $i$.



%

\begin{defn}
An \emph{address} is a string $\sigma \in [N]^\NN$ such that $\sigma_k^+=\sigma_{k+1}^-$ for every $k$. We write $\sigma^-=\sigma_1^-$ and denote the set of addresses by $\Sigma$. For $v\in \V$ we denote by $\Sigma_v$ the set of addresses $\sigma$ with $\sigma^-=v$.
\end{defn}


\begin{rmk}
We can write $\Sigma = \ban\set{ij:i^+\ne j^-}$, immediately showing that $\Sigma$ is a shift space of finite type.
\end{rmk}

A \emph{path} is a word in $\mcL_\Sigma$. For a path $\alpha$ of length $n$, we write $\alpha^-=\alpha_1^-$, $\alpha^+=\alpha_n^+$, and $f_\alpha= f_{\alpha_1}\cdots f_{\alpha_n}$. One should think of the graph $\G$ as determining the orders in which we may compose the maps $f_i$.




\begin{defn}
The \emph{address map} $\pi: \Sigma \to \RR^M$ is defined by
\[\pi(\sigma) = \lim_{k \to \infty} f_{\sigma|_{[1,k]}} (x)\]
for any $x\in \RR^M$. The \emph{attractor} $A$ of the graph IFS is the image $\pi(\Sigma)$ of $\Sigma$ under $\pi$.
\end{defn}

Just as for the classical IFS theory, the address map $\pi:\Sigma\to \RR^M$ is indeed well-defined and continuous.
We may associate to a path $\alpha$ the set $\Sigma_\alpha$ of addresses which begin with $\alpha$. We write $\pi(\alpha)=\pi(\Sigma_\alpha)$.



\begin{note}
We will assume that the components of the attractor are pairwise disjoint, that is, we have $A_v\cap A_{v'}=\varnothing$ whenever $v\ne v'$.
\end{note}

\section{Fractal tops and shifts}

A point in the attractor of an IFS may have multiple addresses. \cite{blowups_tops} gives a way to choose a preferred address for each point. In this section we investigate how this generalises to graph-directed systems, as well as explore some results about the set of top addresses.

\subsection{The fractal top}

We introduce top addresses for graph IFSs and study their connections with shift spaces.

\begin{defn}
We define $\tau:A\to \Sigma$ by $\tau(x)=\min\pi^{-1}(x)$ and we denote by $\Sigma_\tops = \tau(A)$ the set of \emph{top addresses}.
\end{defn}

\begin{rmk}
The above definition is well-defined by Lemma \ref{lem:closed}. The restriction $\pi:\Sigma_\tops\to A$ is a bijection.
\end{rmk}

The proof of the next result is essentially identical to the proof for the classical IFS case, as found in \cite{blowups_tops}.

\begin{lem}\label{lem:contain}
We have $S(\Sigma_\tops) \subseteq \Sigma_\tops$.
\end{lem}

\begin{proof}
Given $\sigma \in \Sigma_\tops$, we have
\begin{align*}
    \sigma &= \tau (\pi (\sigma)) \\
    &= \min \pi^{-1} (\pi (\sigma)) \\
    &= \sigma_1 \min \set{ \omega \in \Sigma : \sigma_1 \omega \in \pi^{-1} (\pi (\sigma)) } \\
    &= \sigma_1 \min \set{ \omega \in \Sigma : \pi(\sigma_1 \omega) = \pi(\sigma) } \\
    &= \sigma_1 \min \set{ \omega \in \Sigma : f_{\sigma_1} (\pi(\omega)) = f_{\sigma_1} (\pi(S(\sigma))) } \\
    &= \sigma_1 \min \set{ \omega \in \Sigma : \pi(\omega) = \pi(S(\sigma)) } \\
    &= \sigma_1 \min \pi^{-1} ( \pi (S(\sigma))) \\
    &= \sigma_1 \tau (\pi (S(\sigma))),
\end{align*}
where the third equality follows by Lemma \ref{lem:minimum} and the sixth equality follows because $f_{\sigma_1}$ is invertible.
Thus
\[ S(\sigma) = S(\sigma_1 \tau(\pi(S(\sigma)))) = \tau (\pi(S(\sigma))) \in \Sigma_\tops.\qedhere \]
\end{proof}

We now come to the main result of this section.

\begin{thm}\label{thm:topshift}
The closure $\overline\Sigma_\tops$ of $\Sigma_\tops$ is a shift space.
\end{thm}

\begin{proof}
Lemma \ref{lem:contain} and the continuity of $S$ imply that $S(\overline\Sigma_\tops)\subseteq \overline\Sigma_\tops$, so the result follows by Proposition \ref{prop:character}.
\end{proof}

This result suggests that it may be of interest to study the closure of the top further. This is the primary concern of the remainder of this section.


\begin{defn}
A collection $\set{U_v}_{v\in \V}$ of nonempty open subsets of $\RR^M$ obeys the \emph{open set condition} for $(\F,\G)$ if
\begin{enumerate}
    \item $f_i(U_{i^+})\subseteq U_{i^-}$ for each $i$
    \item $f_i(U_{i^+})\cap f_j(U_{j^+})=\varnothing$ whenever $i\ne j$.
\end{enumerate}
\end{defn}

\begin{lem}
If $\set{U_v}_{v\in \V}$ obeys the open set condition, then $\overline U_v \supseteq A_v$ for each $v\in \V$.
\end{lem}

\begin{proof}
Let $x\in A_v$ and let $\sigma$ be an address for $x$. We have $\sigma^-=v$. For each $k$, pick $x_k\in f_{\sigma|_{[1,k]}}(U_{(\sigma|_{[1,k]})^+})$. Then $x_k\in U_v$ and $x_k\in \pi(\sigma|_{[1,k]})$ for each $k$. So $x_k\to x$.
\end{proof}

\begin{lem}\label{lem:opensetcond}
Suppose $(\F,\G)$ has a collection $\set{U_v}_{v\in \V}$ obeying the open set condition. Then $\overline\Sigma_\tops = \Sigma$.
\end{lem}

\begin{proof}
Note that $f_i(U_{i^+})\cap f_j(\overline U_{j^+})=\varnothing$ whenever $i\ne j$. It follows by induction that $f_\alpha(U_{\alpha^+})\cap f_\beta(\overline U_{\beta^+})=\varnothing$ whenever $\alpha$ and $\beta$ are distinct paths of equal length.

Now, let $\alpha$ be a path and let $x\in f_\alpha (U_{\alpha^+}\cap A_{\alpha^+})\subseteq f_\alpha(U_{\alpha^+})$. Suppose $\sigma$ is an address of $x$ and write $\sigma=\beta\omega$ where $\abs\alpha=\abs\beta$. Then
\[ x= \pi(\sigma)
=\pi(\beta\omega)= f_\beta\pi(\omega)\in f_\beta(A_{\beta^+})\subseteq f_\beta(\overline U_{\beta^+}). \]
That is, $x\in f_\alpha(U_{\alpha^+})\cap f_\beta(\overline U_{\beta^+})$, so $\alpha=\beta$. In particular, $\tau(\pi(\sigma))$ begins with $\alpha$, so $\alpha$ appears in $\Sigma_\tops$. The result follows by Theorem~\ref{thm:topshift}.
\end{proof}

\begin{defn}
A graph IFS is said to be
\begin{enumerate}
    \item \emph{totally disconnected} if $f_i(A_{i^+})\cap f_j(A_{j^+})=\varnothing$ for every $i$
    \item \emph{just touching} if it is not totally disconnected but has a collection of open sets obeying the open set condition
    \item \emph{overlapping} if $\Int(f_i(A_{i^+})\cap f_j(A_{j^+}))\ne \varnothing$ for some $i\ne j$.
\end{enumerate}
\end{defn}

\begin{thm}\label{thm:classify}
For a graph IFS $(\F,\G)$, we have
\begin{enumerate}
    \item if $(\F,\G)$ is totally disconnected, then $\Sigma_\tops = \Sigma$
    \item if $(\F,\G)$ is just touching, then $\Sigma_\tops \subsetneq \overline\Sigma_\tops = \Sigma$
    \item if $(\F,\G)$ is overlapping, then $\overline\Sigma_\tops \subsetneq \Sigma$.
\end{enumerate}
\end{thm}

\begin{proof}
If $(\F,\G)$ is totally disconnected then $\pi$ is injective, so $\Sigma_\tops = \Sigma$. This proves (a).

If $(\F,\G)$ is just touching, Lemma \ref{lem:opensetcond} gives that $\overline \Sigma_\tops = \Sigma$. But $\pi$ is not injective, so $\Sigma_\tops \subsetneq \overline \Sigma_\tops$, proving (b).

Suppose $(\F,\G)$ is overlapping and let $i<j$ be such that $\Int(f_i(A_{i^+})\cap f_j(A_{j^+}))\ne \varnothing$. Then there exists a path $\alpha$ starting with $j$ such that $\pi(\alpha)\subseteq \Int(f_i(A_{i^+})\cap f_j(A_{j^+}))$. But every point in $\pi(\alpha)$ has an address starting with $i$, so $\alpha$ does not begin any top address. But, by Lemma~\ref{lem:contain}, this implies $\alpha$ does not appear anywhere in $\Sigma_\tops$. But then $\overline\Sigma_\tops$ also does not contain $\alpha$. This proves (c).
\end{proof}

\subsection{Banned words in top shift spaces}

In this subsection we discuss some results about the set of reduced banned words in $\overline{\Sigma}_\tops$.

\begin{prop}\label{prop:redbannedint}
    If $\alpha = \alpha_1\cdots\alpha_n$ is a reduced banned word in $\overline{\Sigma}_\tops$, then $f_\alpha(A) \cap f_k(A) \ne \varnothing$ for some $k < \alpha_1$.
\end{prop}

\begin{proof}
    Let $\alpha = \alpha_1\cdots\alpha_n$ be a reduced banned word. Then $\alpha_2\cdots\alpha_n$ is not banned, so there exists some $x \in A$ with $\tau(x) \in \Sigma_{\alpha_2\cdots\alpha_n}$. Let $\sigma$ be the top address of $f_{\alpha_1}(x)$. Since $\pi(\sigma) = \pi(\alpha_1\tau(x))$, we have that $\sigma_1 \leq \alpha_1$.

    Suppose $\alpha_1 = \sigma_1$. Then $S(\sigma)$ must be an address for $x$. But since $S(\sigma) \in \Sigma_\tops$, we must have that $S(\sigma) = \tau(x)$. So $\sigma|_{[1,n]} = \alpha$, which gives a contradiction. Therefore $\sigma_1 < \alpha_1$. Also note that $f_{\alpha_1}(x) \in f_\alpha(A)\cap f_{\sigma_1}(A)$, so the Proposition holds.
\end{proof}

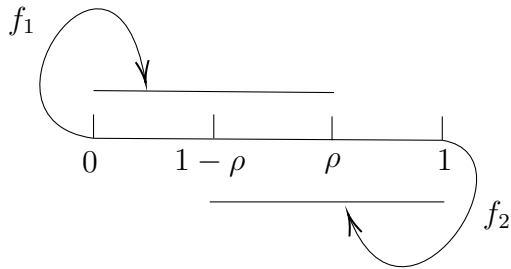
\begin{figure}[h!]
    \centering
    \begin{tikzpicture}[x=0.75pt,y=0.75pt,yscale=-1,xscale=1]

\draw    (88,200) -- (262,201) ;
\draw    (88,188) -- (88,200) ;
\draw    (262,189) -- (262,201) ;
\draw    (148,188) -- (148,200) ;
\draw    (207,189) -- (207,201) ;
\draw    (88,176) -- (208,177) ;
\draw    (146,232) -- (263,232) ;
\draw    (88,200) .. controls (21.33,191.05) and (98.22,79.13) .. (113.77,170.61) ;
\draw [shift={(114,172)}, rotate = 260.93] [color={rgb, 255:red, 0; green, 0; blue, 0 }  ][line width=0.75]    (10.93,-3.29) .. controls (6.95,-1.4) and (3.31,-0.3) .. (0,0) .. controls (3.31,0.3) and (6.95,1.4) .. (10.93,3.29)   ;
\draw    (262,201) .. controls (310.76,208.96) and (242.69,310.97) .. (215.41,239.1) ;
\draw [shift={(215,238)}, rotate = 69.95] [color={rgb, 255:red, 0; green, 0; blue, 0 }  ][line width=0.75]    (10.93,-3.29) .. controls (6.95,-1.4) and (3.31,-0.3) .. (0,0) .. controls (3.31,0.3) and (6.95,1.4) .. (10.93,3.29)   ;

\draw (81,205.4) node [anchor=north west][inner sep=0.75pt]    {$0$};
\draw (258,204.4) node [anchor=north west][inner sep=0.75pt]    {$1$};
\draw (202,206.4) node [anchor=north west][inner sep=0.75pt]    {$\rho $};
\draw (127,204.4) node [anchor=north west][inner sep=0.75pt]    {$1-\rho $};
\draw (44,132.4) node [anchor=north west][inner sep=0.75pt]    {$f_{1}$};
\draw (281,230.4) node [anchor=north west][inner sep=0.75pt]    {$f_{2}$};

\end{tikzpicture}
    \caption{The overlapping IFS with two maps}
    \label{overlap example}
\end{figure}

In general, identifying whether an overlapping GIFS (or even IFS) is of finite type appears to be a difficult problem. For the remainder of this section, we will consider an IFS $\F$ of the form shown in Figure \ref{overlap example} with $ 1/2 \leq \rho < 1$. 

\tikzset{every picture/.style={line width=0.75pt}} 

The system consists of two maps
\begin{equation*}
    f_1(x) = \rho x, \qquad f_2(x) = \rho x + (1 - \rho), 
\end{equation*}

Different choices of $\rho$ will lead to different properties in $\overline{\Sigma}_\tops$.
First we let $\rho = 1/2$. In this case the system is just touching, with the only point in the images of both maps being $1/2$.

\begin{prop}
    When $\rho = 1/2$, we have that $\Sigma_\tops \ne \Sigma$ but $\overline{\Sigma}_\tops = \Sigma$.
\end{prop}

\begin{proof}
    This IFS is just touching, and therefore by Theorem \ref{thm:classify} we have that ${\Sigma_\tops \subsetneq \overline\Sigma_\tops = \Sigma}$.
\end{proof}

Now we will consider the case when $1/2 < \rho < 1$. 

\begin{prop}
    Let $\alpha = \alpha_1 \cdots \alpha_n$ be a word with length $|\alpha| = n$, then
\begin{equation*}
    \begin{split}
        f_\alpha(x) &= \rho^n x + (1 - \rho)[(\alpha_1 - 1) + (\alpha_2 - 1)\rho + \cdots + (\alpha_n - 1) \rho^{n-1}]\\
        & = \rho^n x + (1 - \rho)\sum_{i = 1}^{n} (\alpha_i - 1) \rho^{i-1}. 
    \end{split}
\end{equation*}
Denote $\displaystyle{\alpha(\rho) = \sum_{i = 1}^{n} (\alpha_i - 1) \rho^{i-1}}$. 
\end{prop}

\begin{proof}
    If $n = 1$, then $f_1(x) = \rho x + (1 - \rho)(1 - 1) = \rho x$ and $f_2(x) = \rho x + (1 - \rho)(2 - 1) = \rho x + (1 - \rho)$. Now suppose the formula holds for all words with length less than or equal to $n$. Let $\alpha = \alpha_1\cdots \alpha_{n+1}$ and denote $\beta = \alpha|_{[2, n+1]}$, 
    \begin{equation*}
        f_\beta(x) = \rho^n x + (1 - \rho)\sum_{i = 2}^{n+1} (\alpha_i - 1)\rho^{i-2}. 
    \end{equation*}
    Apply $f_{\alpha_1}$, 
    \begin{equation*}
        \begin{split}
            f_{\alpha}(x) = f_{\alpha_1}\circ f_\beta(x) & = \rho^{n + 1}x + (1 - \rho)\sum_{i = 2}^{n+1} (\alpha_i - 1)\rho^{i-1} + (\alpha_1 - 1)(1 - \rho)\\
            & = \rho^{n + 1}x + (1 - \rho)\sum_{i=1}^{n+1} (\alpha_i - 1)\rho^{i-1}. 
        \end{split}
    \end{equation*}
\end{proof}

We note that $f_\alpha([0,1])$ is merely $[0, \rho^n]$ with a shift $(1-\rho)\alpha(\rho)$. Let $\Sigma^*$ denote the set of all finite words, we have the following characterisation of reduced banned words.

\begin{prop}\label{characterisation of rbw}
    Let $\alpha \in \Sigma^*$ with $|\alpha| = n$, it is a reduced banned word if and only if all of the following hold:
    \begin{enumerate}[label=(A.\arabic*)]
        \item $\alpha_1 = 2$, \label{(a)}

        \item $f_\alpha(1) = \rho^n + (1 - \rho)\alpha(\rho) \leq \rho$, \label{(b)}

        \item For any subword $\beta = \beta_1\cdots \beta_m$ of $\alpha$ with $|\beta| = m$ and $m < n$, we have either $\beta_1 = 1$ or $f_\beta(1) = \rho^{m} + (1 - \rho)\beta(\rho) > \rho$. \label{(c)}
    \end{enumerate}
\end{prop}

\begin{proof}
    Suppose $\alpha$ is reduced banned, \ref{(a)} follows from Proposition \ref{prop:redbannedint}. Note that $\beta$ by definition is not banned, so either $\beta_1 = 1$ or if $\beta_1 = 2$, it must have some points on the top, for which we at least have to require $\rho^{m} + (1 - \rho)\beta(\rho) > \rho$. Hence \ref{(c)} is implied. If $\rho^n + (1 - \rho)\alpha(\rho) > \rho$, then there exists some $x\in [\rho,1]\cap f_\alpha([0,1])$ such that its top address does not start with $\alpha$. But every element in $[\rho, 1]$ has a top address starting with $2$, which indicates that $\alpha|_{[2,n]}$ is not on the top, contradiction. 

    Conversely, \ref{(a)} and \ref{(b)} imply that $\alpha$ is banned. If $\alpha$ satisfies all of \ref{(a)}, \ref{(b)} and \ref{(c)} but is not reduced banned, then there exists a subword $\beta$ of $\alpha$ such that it is a reduced banned word. Previous argument implies that $\beta_1 = 2$ and $\rho^{m} + (1 - \rho)\beta(\rho) \leq \rho$, but this contradicts with \ref{(c)}. 
\end{proof}



Any reduced banned word must end with $1$. 

\begin{lem}\label{rbw ends with 1}
    Let $\alpha$ be any reduced banned word, then $\alpha_{|\alpha|} = 1$. 
\end{lem}

\begin{proof}
    Suppose not, then there exists a reduced banned word $\alpha$ of length $n$ ending with $2$, we note that $1$ is the fixed point of $f_2$. Let $\beta = \alpha|_{[1,n-1]}$, then 
    \begin{equation*}
        f_\alpha(1) = f_\beta(1) \leq \rho
    \end{equation*}
    by \ref{(b)}. This at the same time indicates that $\beta$ is banned, but it contradicts the fact that $\alpha$ is reduced banned.

\end{proof}

Moreover, we find that the IFS always has at least one reduced banned word for any $1/2< \rho< 1$.

\begin{lem}
    The IFS always admits at least one reduced banned word $\alpha = 21\cdots 1$ with $|\alpha| = n$ for some $n > 2$.
\end{lem}

\begin{proof}
    Suppose $\alpha = 21 \cdots 1$ with length $n$, then $\alpha(\rho) = 1$ and 
    \begin{equation*}
        \rho^n + (1 - \rho)\alpha(\rho) - \rho = \rho^n + 1 - 2\rho \leq 0
    \end{equation*}
    for some $n > 2$. 
\end{proof}

Let $m$ denote the length of the word such that $21\cdots 1$ is reduced banned. In fact, it is the first reduced banned word. 

\begin{lem}
    There is no reduced banned word of length less than $m$. 
\end{lem}

\begin{proof}
    By Proposition \ref{characterisation of rbw}, if $\alpha$ is reduced banned, then $f_\alpha([0,1])\subseteq [1 - \rho, \rho]$. Denote the word $21\cdots 1$ of length $n$ by $\Xi_n$, we note that $\Xi_n(\rho) = \min\{\alpha(\rho) : \alpha_1 = 2 , |\alpha| = n\}$. Thus by \ref{(b)}, if $\Xi_n$ is not reduced banned, then $\alpha$ with $|\alpha| = n$ cannot be reduced banned. 
\end{proof}



If $f_{\Xi_m}(1) < \rho$, we observe a property of the second reduced banned word (the reduced banned word of the second smallest length). 

\begin{lem}\label{second reduced banned}
    The second reduced banned word $\beta$ (if exists) starts with $\alpha = \underbrace{21\cdots 12}_m$. 
\end{lem}

\begin{proof}
    We start the search with the smallest possible prefix word $\alpha$, and note that by Lemma \ref{rbw ends with 1} it is not banned. Let $\gamma = \underbrace{21\cdots 1}_{m-1}$, note that $f_\gamma(x) = \rho^{m-1}x + (1 - \rho)$ and 
    \begin{equation*}
        \begin{split}
            f_\gamma^{\circ n}(x) &= \rho^{n(m-1)}x + (1 - \rho)(1 + \rho^{m-1} + \cdots + \rho^{(n-1)(m-1)})\\
            & = \rho^{n(m-1)}x + (1 - \rho)\sum_{i=0}^{n-1} \rho^{i(m-1)}\\
            & = \rho^{n(m-1)}x + (1- \rho) \frac{1 - \rho^{n(m-1)}}{1 - \rho^{m-1}}. 
        \end{split}
    \end{equation*}
    Concatenating $\gamma$ with itself will give the smallest non-banned word that we can find (since the concatenation cannot contain enough consecutive $1$'s to coincide with $\Xi_m$), note that 
    \begin{equation*}
        \lim_{n\to \infty} f_\gamma^{\circ n}(1) = \frac{1 - \rho}{1 - \rho^{m-1}}. 
    \end{equation*}
    Since $\Xi_m$ is reduced banned and by assumption $f_{\Xi_m}(1) < \rho$, thus $\rho^m + (1 - \rho) < \rho$, i.e., $1 - \rho < \rho(1 - \rho^{m-1})$. Hence, 
    \begin{equation*}
        \lim_{n\to \infty} f_\gamma^{\circ n}(1) = \frac{1 - \rho}{1 - \rho^{m-1}} < \rho. 
    \end{equation*}
    Since $f_\gamma$ is a contraction, there exists some $n\in \NN$ such that $f_\gamma^{\circ n}(1) \leq \rho$, which indicates that there is a reduced banned word $\beta$ starting with $\alpha$. 
\end{proof}

In particular, the second reduced banned word $\beta$ (if exists) is found to be of the form $\underbrace{\gamma\cdots \gamma}_{n}\gamma|_{[1, k]}$ for some $n\geq 1$ and $1 < k < m - 1$. 
We then observe that the endpoint $f_\beta(1)$ is to the right of $f_{\Xi_m}(1)$. 

\begin{lem}\label{next rdb is to the right}
    $f_\beta(1) > f_{\Xi_m}(1)$.
\end{lem}

\begin{proof}
    Consider the tail $\gamma = \beta|_{[m, |\beta|]}$, it is not banned, hence $f_\gamma(1) > \rho$. Let $\delta = \Xi_{m-1}$, we note that $f_{\Xi_m}(1) = f_\delta(\rho)$ and $f_\delta$ is monotone. 
    \begin{equation*}
        f_\beta(1) = f_\delta \circ f_\gamma(1) > f_\delta(\rho) = f_{\Xi_m}(1). 
    \end{equation*}
\end{proof}

It indicates that the second reduced banned word does \emph{not} exist if $f_{\Xi_m}(1) = \rho$.

\begin{ex}
    When $\rho = \frac{\sqrt{5} - 1}{2}$, we have $\rho^2 + \rho = 1$ and there is only one reduced banned word $211$.


\end{ex}

Moreover, we know $\beta$ ends with $1$ by Lemma \ref{rbw ends with 1}. If the third reduced banned word $\delta$ exists, then it must start with $\beta|_{[1, |\beta|-1]}2$ as it is the only possible prefix. Likewise, by the exact same argument in Lemma \ref{next rdb is to the right}, we can conclude that $f_{\delta}(1) > f_{\beta}(1)$.

A similar argument can be applied to any reduced banned word in this IFS, hence we deduce 

\begin{cor}
    There are no reduced banned words of the same length. 
\end{cor}

\begin{cor}
    If there is a reduced banned word $\alpha$ such that \ref{(b)} is satisfied with equality, then the IFS cannot have any reduced banned word of length greater than $|\alpha|$. 
\end{cor}


In particular, it gives a sufficient condition for this IFS to be of finite type. 

\begin{cor}\label{sufficient condition}
    The IFS is of finite type if there exists a reduced banned word $\alpha$ such that \ref{(b)} is satisfied with equality.
\end{cor}



\begin{rmk}
    Corollary \ref{sufficient condition}, Proposition \ref{characterisation of rbw} and Lemma \ref{rbw ends with 1} indicate that the IFS is of finite type if the scaling factor $\rho$ is a solution to a monic polynomial of the form 
    \begin{equation*}
        \rho^n + (1 - \rho)\alpha(\rho) = \rho
    \end{equation*}
    for some reduced banned word $\alpha$ of length $n$, where the polynomial has coefficients in $\{0, \pm 1\}$. 
\end{rmk}

We suspect that the converse of Corollary \ref{sufficient condition} is also true, as motivated by the following example. 

\begin{ex}
    When $\rho = 2/3$, the IFS does not satisfy the condition in Corollary \ref{sufficient condition} and we suspect it is of infinite type. A computer program was used to find reduced banned words, up to length $11$, they are
\begin{align*}
    &211\\
    &212121\\
    &2121221\\
    &21212221\\
    &212122221\\
    &212122222121\\
    &212122222122121.
\end{align*}
Let $\{\alpha^i\}_{i\in I}$ be an enumeration of the reduced banned words ordered from the smallest length to the largest length, where $I\subseteq \NN$. Let $\gamma^i = \alpha^i|_{[1,|\alpha^i|-1]}$. From the pattern of these reduced words, it seems that for each $i > 1$, $\alpha^i = \underbrace{\gamma^{i-1}\cdots\gamma^{i-1}}_j\gamma^{i-1}|_{[1, k]}$ for some $j\geq 1$, $1 < k < |\gamma^{i-1}| - 1$ or $\alpha^i = \underbrace{\gamma^{i-1}\cdots\gamma^{i-1}}_j$ for some $j \geq 1$. 
\end{ex}

We conjecture that the pattern holds for all scaling factors $1/2 < \rho < 1$, for which we need the following
\begin{conjecture}
    For any $i\in I$, $\gamma^i\gamma^i$ does not contain any $\alpha^j$ for any $j\leq i$ as a subword.
\end{conjecture}

It turns out to be very difficult to prove or find a counterexample. But if it were true, we can follow a similar argument as in the proof of Lemma \ref{second reduced banned} to show the converse of Corollary \ref{sufficient condition}, hence giving a classification of the scaling factor $\rho$ when this particular IFS is of finite type.

\section{Shift invariance of the top}

In \cite{blowups_tops}, it is proved that $\Sigma_\tops$ is shift invariant for any classical IFS. Our Theorem \ref{thm:topshift} demonstrated one of the inclusions for any graph IFS, but in this section we show that the opposite inclusion is, in general, false in this setting. We explore the behaviours of $\Sigma_\tops$ with respect to the shift operator in several examples.


\subsection{A range of behaviours}\label{diagram section}

This subsection is devoted to examples. We begin with a simple example of a graph IFS for which $\Sigma_\tops$ is not shift invariant.

\begin{ex}\label{ex:not_si}
Let $(\F,\G)$ be the graph IFS on $\RR$ shown in Figure \ref{fig:notinvariant}.
We have $A_1=[0,1]$ and $A_2=[2,3]$, and the system consists of the four maps
\begin{align*}
    f_1(x) &= x/2\\
    f_2(x) &= x/2-1/2\\
    f_3(x) &= x/2+2\\
    f_4(x) &= x/2 + 1.
\end{align*}
We note that $\tau(2) = 3\overline 1$ so that $3\overline 1\in \Sigma_\tops$. But neither $13\overline 1$ nor $33\overline1$ are paths, and $\pi(23\overline 1)=1/2=\pi(12\overline 4)$ and $\pi(43\overline1)=5/2=\pi(32\overline4)$, so $12\overline 1,32\overline4\notin \Sigma_\tops$ because $12\overline 4 < 23\overline 1$ and $32\overline4<43\overline1$. So $3\overline1\notin S(\Sigma_\tops)$.
\end{ex}

%

%


As it turns out, the labels for Example \ref{ex:not_si} can be rearranged to obtain a system in which $\Sigma_\tops$ is indeed shift invariant, as demonstrated in the next example.

\begin{ex}
Let $(\F,\G)$ be the graph IFS shown in Figure \ref{fig:invariant}.
To see that $\Sigma_\tops$ is fact shift invariant in this case, observe that if $\sigma$ is a top address of a point in $\cc01$, then $1\sigma$ is a top address. Similarly, if $\sigma$ is a top address of a point in $\cc23$, then $2\sigma$ is also a top address. So $S(\Sigma_\tops)=\Sigma_\tops$.
\end{ex}









\begin{rmk}
Suppose a graph IFS is such that for each $v\in\V$, there is an edge $i$ with $i^-=i^+=v$. By ordering such that these edges come first, we obtain a system in which $\Sigma_\tops$ is shift invariant.
\end{rmk}

Next, we see a nontrivial example in which $\Sigma_\tops$ is in fact shift invariant for every ordering of the edges.

\begin{ex}\label{ex:ssi}
Let $(\F,\G)$ be the graph IFS shown in Figure \ref{fig:ssi}. Regardless of the ordering of the edges, we can concatenate $i$ with a top address for $0$, we can concatenate $l$ with a top address for $1$, we can concatenate $j$ with a top address for $2$ and we can concatenate $k$ with a top address for $3$. In each case we must obtain a top address. It is easy to see that top addresses of points in the interiors have preimages in the top, so $\Sigma_\tops$ is shift invariant.
\end{ex}

The following is another such example, this time with no edges $i$ with $i^-=i^+$.

\begin{ex}
Let $(\F,\G)$ be the graph IFS shown in Figure \ref{fig:ssi no self ref}. A similar argument to that for Example \ref{ex:ssi} shows that $\Sigma_\tops$ is again shift invariant for any ordering of the edges.
\end{ex}


Finally, we realise a pathological case in which $\Sigma_\tops$ is not shift invariant for any ordering of the edges.

\begin{ex}
Let $(\F,\G)$ be the graph IFS shown in Figure \ref{fig:neverinvariant}.
For the given labelling, it is not possible for both $27\overline6$ and $13\overline4$ to be in $\Sigma_\tops$. These are the only possible preimages of $7\overline 6$ and $3\overline4$, so $\Sigma_\tops$ cannot be shift invariant for any ordering of the edges.
\end{ex}


\subsection{Points not in the shift}

The missing top addresses from the shift map $S$ is the collection $\Upsilon_1 = \{\s \in \Sigma_\tops: \s\not\in S(\Sigma_\tops)\}$. We have the following characterisation. 

\begin{prop}\label{prop:up1}
    $\Upsilon_1$ is equal to the set 
    \begin{align*}
        \set[\bigg]{\s\in \Sigma_\tops: \pi(\s)\in \bigcup_{v\in \V} \bigcap_{j: j^+ = v} f_j^{-1}(f_j(A_v)\cap\bigcup_{k: k^- = j^-, \atop k< j } f_k(A_{k^+}))}.
    \end{align*}
\end{prop}

\begin{proof}
First, let $\sigma \in \Sigma_\tops$ such that $\sigma \notin S(\Sigma_\tops)$. Suppose $\sigma^- = v$. Let $j$ be an edge with $j^+ = v$. Then $j\sigma \notin \Sigma_\tops$, so we have $\pi(j\sigma) = \pi(\omega)$ for some $\omega \in \Sigma$ with $\omega < \sigma$. If $\omega_1 = j$ then $\pi(\sigma) = \pi(\omega_2\cdots)$. Since $\sigma$ is in $\Sigma_\tops$, we must have that $\sigma \leq \omega_2\cdots$. But then $j\sigma \leq \omega$, which gives a contradiction. Therefore $\omega_1 < j$. So $f_j(\pi(\sigma))$ is an element of $f_j(A_v)\cap\bigcup_{k: k^- = j^-, \atop k< j } f_k(A_{k^+})$. This shows the inclusion of $\Upsilon_1$ into the proposed set.

Now let $\sigma$ be an element of the proposed set with $\pi(\sigma) \in A_v$. Let $j$ be an edge with $j^+ = v$. Then $\pi(\sigma) \in f_j^{-1}(f_j(A_v)\cap\bigcup_{k: k^- = j^-, \atop k< j } f_k(A_{k^+}))$. So $\pi(j\sigma) \in \bigcup_{k: k^- = j^-, \atop k< j } f_k(A_{k^+}))$. Therefore there exists some edge $k$ with $k < j$ and $k^-  = j^-$ such that $\pi(j\sigma) \in f_k(A_{k^+})$. So there exists an address $k\omega$ with $\pi(k\omega) = \pi(j\sigma)$ and $k\omega < j\sigma$. So $j\sigma \notin \Sigma_\tops$ and therefore $\sigma \notin S(\Sigma_\tops)$.
\end{proof}




Analogously, we define $\Upsilon_n = \{\s\in \Sigma_\tops: \s\not\in S^n(\Sigma_\tops)\}$ to be the top addresses that do not appear after applying the shift map $n$ times. It is clear that $\Upsilon_n \subseteq \Upsilon_{n+1}$. Moreover, 





\begin{prop}
    $\Upsilon_{n+1}$ is equal to the set 
    \begin{equation*}
        \Upsilon_n \cup \set[\bigg]{\s\in \Sigma_\tops: \pi(\s)\in \bigcup_{v\in \V} \bigcap_{\alpha: \alpha^+ = v, \atop |\alpha| = n+1} f_{\alpha}^{-1}(f_{\alpha}(A_v)\cap \bigcup_{k: k^- = \alpha^-, \atop k < \alpha_1} f_{k}(A_{k^+}))}.
    \end{equation*}
\end{prop}

This can be shown using an argument similar to that in the proof of Proposition \ref{prop:up1}.

\section*{Acknowledgements}
This research was majorly done while we were doing a mathematics Honours at ANU, we would like to thank Louisa F. Barnsley and Michael F. Barnsley for offering a special topics course on fractal tiling theory, as well as their invaluable guidance throughout the completion of this paper.

\clearpage

\appendix

\section{Diagrams}

We collect the diagrams needed for Section \ref{diagram section}.
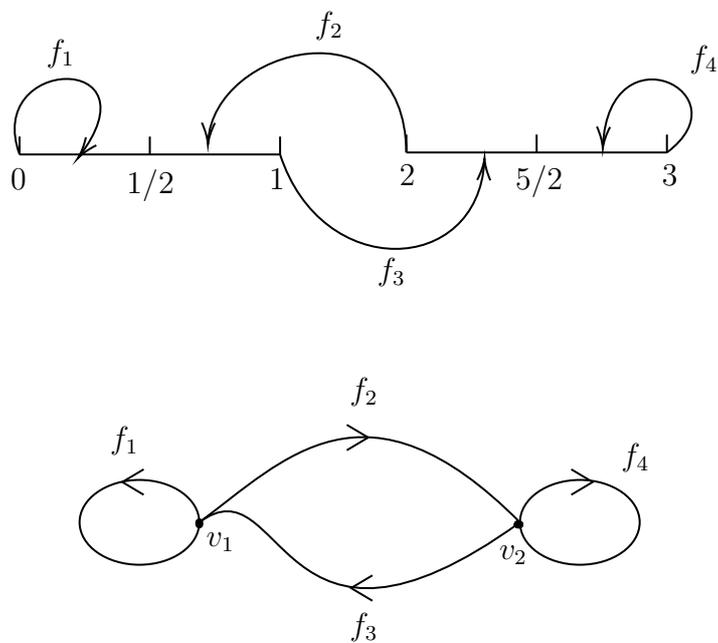
\begin{figure}[h!]
    \centering
    \begin{subfigure}{1\textwidth}
        \centering
        \begin{tikzpicture}[x=0.75pt,y=0.75pt,yscale=-1,xscale=1]

\draw    (61,229) -- (191,229) ;
\draw    (126,220) -- (126,229) ;
\draw    (61,220) -- (61,229) ;
\draw    (191,220) -- (191,229) ;
\draw    (254,228) -- (384,228) ;
\draw    (319,219) -- (319,228) ;
\draw    (254,219) -- (254,228) ;
\draw    (384,219) -- (384,228) ;
\draw    (61,229) .. controls (43.09,179.25) and (131.11,177.02) .. (90.62,230.19) ;
\draw [shift={(90,231)}, rotate = 307.87] [color={rgb, 255:red, 0; green, 0; blue, 0 }  ][line width=0.75]    (10.93,-3.29) .. controls (6.95,-1.4) and (3.31,-0.3) .. (0,0) .. controls (3.31,0.3) and (6.95,1.4) .. (10.93,3.29)   ;
\draw    (254,228) .. controls (253.01,150.78) and (157.93,174.51) .. (155.06,222.54) ;
\draw [shift={(155,224)}, rotate = 271.17] [color={rgb, 255:red, 0; green, 0; blue, 0 }  ][line width=0.75]    (10.93,-3.29) .. controls (6.95,-1.4) and (3.31,-0.3) .. (0,0) .. controls (3.31,0.3) and (6.95,1.4) .. (10.93,3.29)   ;
\draw    (384,228) .. controls (423.6,198.3) and (355.39,164.68) .. (352.08,224.17) ;
\draw [shift={(352,226)}, rotate = 271.85] [color={rgb, 255:red, 0; green, 0; blue, 0 }  ][line width=0.75]    (10.93,-3.29) .. controls (6.95,-1.4) and (3.31,-0.3) .. (0,0) .. controls (3.31,0.3) and (6.95,1.4) .. (10.93,3.29)   ;
\draw    (191,229) .. controls (211.79,293.35) and (289.43,290.07) .. (292.92,232.75) ;
\draw [shift={(293,231)}, rotate = 91.94] [color={rgb, 255:red, 0; green, 0; blue, 0 }  ][line width=0.75]    (10.93,-3.29) .. controls (6.95,-1.4) and (3.31,-0.3) .. (0,0) .. controls (3.31,0.3) and (6.95,1.4) .. (10.93,3.29)   ;

\draw (113,235) node [anchor=north west][inner sep=0.75pt]   [align=left] {1/2};
\draw (55,235) node [anchor=north west][inner sep=0.75pt]   [align=left] {0};
\draw (184,235) node [anchor=north west][inner sep=0.75pt]   [align=left] {1};
\draw (249,233) node [anchor=north west][inner sep=0.75pt]   [align=left] {2};
\draw (307,233) node [anchor=north west][inner sep=0.75pt]   [align=left] {5/2};
\draw (380,233) node [anchor=north west][inner sep=0.75pt]   [align=left] {3};
\draw (73,169.4) node [anchor=north west][inner sep=0.75pt]    {$f_{1}$};
\draw (207,155.4) node [anchor=north west][inner sep=0.75pt]    {$f_{2}$};
\draw (394,172.4) node [anchor=north west][inner sep=0.75pt]    {$f_{4}$};
\draw (238,279.4) node [anchor=north west][inner sep=0.75pt]    {$f_{3}$};

\end{tikzpicture}
    \end{subfigure}
    \hfill
    \vspace{0.5cm}
    \begin{subfigure}{1\textwidth}       
    \centering
    \begin{tikzpicture}[x=0.8pt,y=0.8pt,yscale=-1,xscale=1]

\draw  [line width=3] [line join = round][line cap = round] (111,150) .. controls (111,150.33) and (111,150.67) .. (111,151) ;
\draw  [line width=3] [line join = round][line cap = round] (260,151) .. controls (261.33,151) and (261.33,151) .. (260,151) ;
\draw   (180,104) -- (190,110) -- (180,116) ;
\draw    (111,150) .. controls (151,120) and (189,77) .. (261,150) ;
\draw    (111,150) .. controls (151,120) and (147,233) .. (261,150) ;
\draw   (192,175) -- (182,181) -- (192,187) ;
\draw   (55,150) .. controls (55,138.95) and (67.54,130) .. (83,130) .. controls (98.46,130) and (111,138.95) .. (111,150) .. controls (111,161.05) and (98.46,170) .. (83,170) .. controls (67.54,170) and (55,161.05) .. (55,150) -- cycle ;
\draw   (85,125) -- (75,131) -- (85,137) ;
\draw   (261,150) .. controls (261,138.95) and (273.54,130) .. (289,130) .. controls (304.46,130) and (317,138.95) .. (317,150) .. controls (317,161.05) and (304.46,170) .. (289,170) .. controls (273.54,170) and (261,161.05) .. (261,150) -- cycle ;
\draw   (285,125) -- (295,131) -- (285,137) ;

\draw (113,153.4) node [anchor=north west][inner sep=0.75pt]    {$v_{1}$};
\draw (250,159.4) node [anchor=north west][inner sep=0.75pt]    {$v_{2}$};
\draw (68,103.4) node [anchor=north west][inner sep=0.75pt]    {$f_{1}$};
\draw (180,80.4) node [anchor=north west][inner sep=0.75pt]    {$f_{2}$};
\draw (180,191.4) node [anchor=north west][inner sep=0.75pt]    {$f_{3}$};
\draw (307,111.4) node [anchor=north west][inner sep=0.75pt]    {$f_{4}$};

\end{tikzpicture}
    \end{subfigure}
    \caption{A graph IFS such that $\Sigma_\tops$ is not shift invariant}
    \label{fig:notinvariant}
\end{figure}

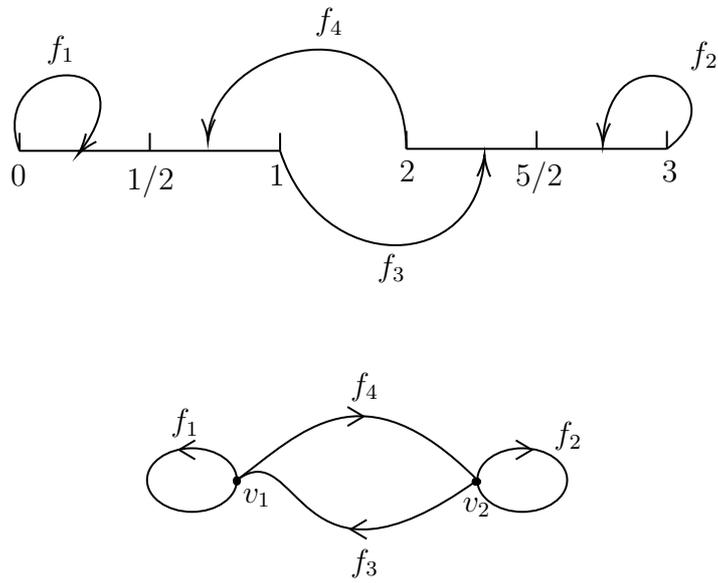
\begin{figure}
\centering
    \begin{subfigure}{1\textwidth}
    \centering
    \begin{tikzpicture}[x=0.75pt,y=0.75pt,yscale=-1,xscale=1]

    \draw    (61,229) -- (191,229) ;
    \draw    (126,220) -- (126,229) ;
    \draw    (61,220) -- (61,229) ;
    \draw    (191,220) -- (191,229) ;
\draw    (254,228) -- (384,228) ;
\draw    (319,219) -- (319,228) ;
\draw    (254,219) -- (254,228) ;
\draw    (384,219) -- (384,228) ;
\draw    (61,229) .. controls (43.09,179.25) and (131.11,177.02) .. (90.62,230.19) ;
\draw [shift={(90,231)}, rotate = 307.87] [color={rgb, 255:red, 0; green, 0; blue, 0 }  ][line width=0.75]    (10.93,-3.29) .. controls (6.95,-1.4) and (3.31,-0.3) .. (0,0) .. controls (3.31,0.3) and (6.95,1.4) .. (10.93,3.29)   ;
\draw    (254,228) .. controls (253.01,150.78) and (157.93,174.51) .. (155.06,222.54) ;
\draw [shift={(155,224)}, rotate = 271.17] [color={rgb, 255:red, 0; green, 0; blue, 0 }  ][line width=0.75]    (10.93,-3.29) .. controls (6.95,-1.4) and (3.31,-0.3) .. (0,0) .. controls (3.31,0.3) and (6.95,1.4) .. (10.93,3.29)   ;
\draw    (384,228) .. controls (423.6,198.3) and (355.39,164.68) .. (352.08,224.17) ;
\draw [shift={(352,226)}, rotate = 271.85] [color={rgb, 255:red, 0; green, 0; blue, 0 }  ][line width=0.75]    (10.93,-3.29) .. controls (6.95,-1.4) and (3.31,-0.3) .. (0,0) .. controls (3.31,0.3) and (6.95,1.4) .. (10.93,3.29)   ;
\draw    (191,229) .. controls (211.79,293.35) and (289.43,290.07) .. (292.92,232.75) ;
\draw [shift={(293,231)}, rotate = 91.94] [color={rgb, 255:red, 0; green, 0; blue, 0 }  ][line width=0.75]    (10.93,-3.29) .. controls (6.95,-1.4) and (3.31,-0.3) .. (0,0) .. controls (3.31,0.3) and (6.95,1.4) .. (10.93,3.29)   ;

\draw (113,235) node [anchor=north west][inner sep=0.75pt]   [align=left] {1/2};
\draw (55,235) node [anchor=north west][inner sep=0.75pt]   [align=left] {0};
\draw (184,235) node [anchor=north west][inner sep=0.75pt]   [align=left] {1};
\draw (249,233) node [anchor=north west][inner sep=0.75pt]   [align=left] {2};
\draw (307,233) node [anchor=north west][inner sep=0.75pt]   [align=left] {5/2};
\draw (380,233) node [anchor=north west][inner sep=0.75pt]   [align=left] {3};
\draw (73,169.4) node [anchor=north west][inner sep=0.75pt]    {$f_{1}$};
\draw (207,155.4) node [anchor=north west][inner sep=0.75pt]    {$f_{4}$};
\draw (394,172.4) node [anchor=north west][inner sep=0.75pt]    {$f_{2}$};
\draw (238,279.4) node [anchor=north west][inner sep=0.75pt]    {$f_{3}$};

\end{tikzpicture}
\end{subfigure}
    \hfill
    \vspace{0.5cm}
    \begin{subfigure}{1\textwidth} 
\centering
\begin{tikzpicture}[x=0.6pt,y=0.6pt,yscale=-1,xscale=1]

\draw  [line width=3] [line join = round][line cap = round] (111,150) .. controls (111,150.33) and (111,150.67) .. (111,151) ;
\draw  [line width=3] [line join = round][line cap = round] (260,151) .. controls (261.33,151) and (261.33,151) .. (260,151) ;
\draw   (180,104) -- (190,110) -- (180,116) ;
\draw    (111,150) .. controls (151,120) and (189,77) .. (261,150) ;
\draw    (111,150) .. controls (151,120) and (147,233) .. (261,150) ;
\draw   (192,175) -- (182,181) -- (192,187) ;
\draw   (55,150) .. controls (55,138.95) and (67.54,130) .. (83,130) .. controls (98.46,130) and (111,138.95) .. (111,150) .. controls (111,161.05) and (98.46,170) .. (83,170) .. controls (67.54,170) and (55,161.05) .. (55,150) -- cycle ;
\draw   (85,125) -- (75,131) -- (85,137) ;
\draw   (261,150) .. controls (261,138.95) and (273.54,130) .. (289,130) .. controls (304.46,130) and (317,138.95) .. (317,150) .. controls (317,161.05) and (304.46,170) .. (289,170) .. controls (273.54,170) and (261,161.05) .. (261,150) -- cycle ;
\draw   (285,125) -- (295,131) -- (285,137) ;

\draw (113,153.4) node [anchor=north west][inner sep=0.75pt]    {$v_{1}$};
\draw (250,159.4) node [anchor=north west][inner sep=0.75pt]    {$v_{2}$};
\draw (68,103.4) node [anchor=north west][inner sep=0.75pt]    {$f_{1}$};
\draw (180,80.4) node [anchor=north west][inner sep=0.75pt]    {$f_{4}$};
\draw (180,191.4) node [anchor=north west][inner sep=0.75pt]    {$f_{3}$};
\draw (307,111.4) node [anchor=north west][inner sep=0.75pt]    {$f_{2}$};

\end{tikzpicture}
\end{subfigure}

\caption{A relabelling of the graph IFS in figure \ref{fig:notinvariant} such that $\Sigma_\tops$ is shift invariant}
\label{fig:invariant}
\end{figure}

\begin{figure}
    \centering
    
    \begin{subfigure}{1\textwidth}
        \centering
        \tikzset{every picture/.style={line width=0.75pt}} 

        \begin{tikzpicture}[x=0.75pt,y=0.75pt,yscale=-1,xscale=1]

\draw    (52,189) -- (185,190) ;
\draw    (52,177) -- (52,189) ;
\draw    (185,178) -- (185,190) ;
\draw    (52,189) .. controls (-14.67,180.05) and (70.14,93.87) .. (85.77,185.6) ;
\draw [shift={(86,187)}, rotate = 260.93] [color={rgb, 255:red, 0; green, 0; blue, 0 }  ][line width=0.75]    (10.93,-3.29) .. controls (6.95,-1.4) and (3.31,-0.3) .. (0,0) .. controls (3.31,0.3) and (6.95,1.4) .. (10.93,3.29)   ;
\draw    (276,190) .. controls (228.24,133.28) and (160.68,103.3) .. (153.11,178.85) ;
\draw [shift={(153,180)}, rotate = 275.19] [color={rgb, 255:red, 0; green, 0; blue, 0 }  ][line width=0.75]    (10.93,-3.29) .. controls (6.95,-1.4) and (3.31,-0.3) .. (0,0) .. controls (3.31,0.3) and (6.95,1.4) .. (10.93,3.29)   ;
\draw    (276,190) -- (409,191) ;
\draw    (276,178) -- (276,190) ;
\draw    (409,179) -- (409,191) ;
\draw    (276,190) .. controls (244.16,141.24) and (314.29,94.47) .. (308.1,185.62) ;
\draw [shift={(308,187)}, rotate = 274.3] [color={rgb, 255:red, 0; green, 0; blue, 0 }  ][line width=0.75]    (10.93,-3.29) .. controls (6.95,-1.4) and (3.31,-0.3) .. (0,0) .. controls (3.31,0.3) and (6.95,1.4) .. (10.93,3.29)   ;
\draw    (118.5,177.5) -- (118.5,189.5) ;
\draw    (342.5,178.5) -- (342.5,190.5) ;
\draw    (185,190) .. controls (250.67,292.49) and (335.15,243.5) .. (381.31,191.78) ;
\draw [shift={(382,191)}, rotate = 131.5] [color={rgb, 255:red, 0; green, 0; blue, 0 }  ][line width=0.75]    (10.93,-3.29) .. controls (6.95,-1.4) and (3.31,-0.3) .. (0,0) .. controls (3.31,0.3) and (6.95,1.4) .. (10.93,3.29)   ;

\draw (45,194.4) node [anchor=north west][inner sep=0.75pt]    {$0$};
\draw (180,195.4) node [anchor=north west][inner sep=0.75pt]    {$1$};
\draw (8,121.4) node [anchor=north west][inner sep=0.75pt]    {$f_{i}$};
\draw (173,109.4) node [anchor=north west][inner sep=0.75pt]    {$f_{k}$};
\draw (269,195.4) node [anchor=north west][inner sep=0.75pt]    {$2$};
\draw (403,195.4) node [anchor=north west][inner sep=0.75pt]    {$3$};
\draw (281,108.4) node [anchor=north west][inner sep=0.75pt]    {$f_{j}$};
\draw (335,238.4) node [anchor=north west][inner sep=0.75pt]    {$f_{l}$};
\end{tikzpicture}
    \end{subfigure}
    
    \hfill
    \vspace{0.5cm}
    \begin{subfigure}{1\textwidth}
        \centering
        \tikzset{every picture/.style={line width=0.75pt}} 
        \begin{tikzpicture}[x=0.6pt,y=0.6pt,yscale=-1,xscale=1]

\draw  [line width=3] [line join = round][line cap = round] (111,150) .. controls (111,150.33) and (111,150.67) .. (111,151) ;
\draw  [line width=3] [line join = round][line cap = round] (260,151) .. controls (261.33,151) and (261.33,151) .. (260,151) ;
\draw   (180,104) -- (190,110) -- (180,116) ;
\draw    (111,150) .. controls (151,120) and (189,77) .. (261,150) ;
\draw    (111,150) .. controls (151,120) and (147,233) .. (261,150) ;
\draw   (192,175) -- (182,181) -- (192,187) ;
\draw   (55,150) .. controls (55,138.95) and (67.54,130) .. (83,130) .. controls (98.46,130) and (111,138.95) .. (111,150) .. controls (111,161.05) and (98.46,170) .. (83,170) .. controls (67.54,170) and (55,161.05) .. (55,150) -- cycle ;
\draw   (85,125) -- (75,131) -- (85,137) ;
\draw   (261,150) .. controls (261,138.95) and (273.54,130) .. (289,130) .. controls (304.46,130) and (317,138.95) .. (317,150) .. controls (317,161.05) and (304.46,170) .. (289,170) .. controls (273.54,170) and (261,161.05) .. (261,150) -- cycle ;
\draw   (285,125) -- (295,131) -- (285,137) ;

\draw (113,153.4) node [anchor=north west][inner sep=0.75pt]    {$v_{1}$};
\draw (250,159.4) node [anchor=north west][inner sep=0.75pt]    {$v_{2}$};
\draw (68,103.4) node [anchor=north west][inner sep=0.75pt]    {$f_{i}$};
\draw (180,80.4) node [anchor=north west][inner sep=0.75pt]    {$f_{k}$};
\draw (180,191.4) node [anchor=north west][inner sep=0.75pt]    {$f_{l}$};
\draw (307,111.4) node [anchor=north west][inner sep=0.75pt]    {$f_{j}$};

\end{tikzpicture}
    \end{subfigure}
    
    \caption{A graph IFS for which $\Sigma_\tops$ is shift invariant for every ordering of the edges}
    \label{fig:ssi}
\end{figure}

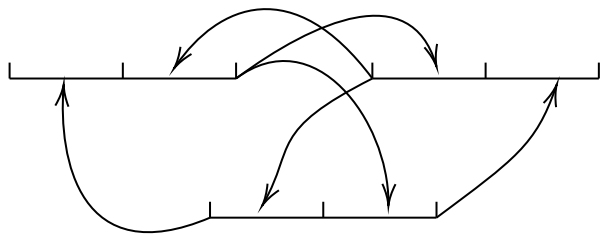
\begin{figure}
        \centering

\tikzset{every picture/.style={line width=0.75pt}} 

\begin{tikzpicture}[x=0.75pt,y=0.75pt,yscale=-1,xscale=1]

\draw    (99,120) -- (212,120) ;
\draw    (280,120) -- (393,120) ;
\draw    (199,190) -- (312,190) ;
\draw    (99,120) -- (99,112) ;
\draw    (212,120) -- (212,112) ;
\draw    (280,120) -- (280,112) ;
\draw    (393,120) -- (393,112) ;
\draw    (199,190) -- (199,182) ;
\draw    (312,190) -- (312,182) ;
\draw    (155.5,120) -- (155.5,112) ;
\draw    (336.5,120) -- (336.5,112) ;
\draw    (255.5,190) -- (255.5,182) ;
\draw    (199,190) .. controls (133.99,217.58) and (123.31,160.75) .. (125.87,124.64) ;
\draw [shift={(126,123)}, rotate = 94.76] [color={rgb, 255:red, 0; green, 0; blue, 0 }  ][line width=0.75]    (10.93,-3.29) .. controls (6.95,-1.4) and (3.31,-0.3) .. (0,0) .. controls (3.31,0.3) and (6.95,1.4) .. (10.93,3.29)   ;
\draw    (312,190) .. controls (351.2,160.6) and (361.59,153.29) .. (371.4,124.77) ;
\draw [shift={(372,123)}, rotate = 108.43] [color={rgb, 255:red, 0; green, 0; blue, 0 }  ][line width=0.75]    (10.93,-3.29) .. controls (6.95,-1.4) and (3.31,-0.3) .. (0,0) .. controls (3.31,0.3) and (6.95,1.4) .. (10.93,3.29)   ;
\draw    (212,120) .. controls (251.6,90.3) and (299.04,70.4) .. (311.63,112.7) ;
\draw [shift={(312,114)}, rotate = 254.74] [color={rgb, 255:red, 0; green, 0; blue, 0 }  ][line width=0.75]    (10.93,-3.29) .. controls (6.95,-1.4) and (3.31,-0.3) .. (0,0) .. controls (3.31,0.3) and (6.95,1.4) .. (10.93,3.29)   ;
\draw    (280,120) .. controls (240.6,66.81) and (203.14,82.51) .. (181.96,113.57) ;
\draw [shift={(181,115)}, rotate = 303.27] [color={rgb, 255:red, 0; green, 0; blue, 0 }  ][line width=0.75]    (10.93,-3.29) .. controls (6.95,-1.4) and (3.31,-0.3) .. (0,0) .. controls (3.31,0.3) and (6.95,1.4) .. (10.93,3.29)   ;
\draw    (212,120) .. controls (251.4,90.45) and (287.89,141.43) .. (288.02,182.15) ;
\draw [shift={(288,184)}, rotate = 271.4] [color={rgb, 255:red, 0; green, 0; blue, 0 }  ][line width=0.75]    (10.93,-3.29) .. controls (6.95,-1.4) and (3.31,-0.3) .. (0,0) .. controls (3.31,0.3) and (6.95,1.4) .. (10.93,3.29)   ;
\draw    (280,120) .. controls (228.06,146.46) and (241.43,154.67) .. (225.98,182.28) ;
\draw [shift={(225,184)}, rotate = 300.38] [color={rgb, 255:red, 0; green, 0; blue, 0 }  ][line width=0.75]    (10.93,-3.29) .. controls (6.95,-1.4) and (3.31,-0.3) .. (0,0) .. controls (3.31,0.3) and (6.95,1.4) .. (10.93,3.29)   ;

\end{tikzpicture}

        \caption{A graph IFS with no self-referencing maps for which $\Sigma_\tops$ is shift invariant for every ordering}
        \label{fig:ssi no self ref}
    \end{figure}

\input{nsi}


\begin{thebibliography}{1}

\bibitem{blowups_tops}
Louisa~F. Barnsley and Michael~F. Barnsley.
\newblock Blowups and tops of overlapping iterated function systems, 2023.

\bibitem{fractals_everywhere}
Michael~F. Barnsley.
\newblock {\em Fractals Everywhere}.
\newblock Academic Press, second edition, 1993.

\bibitem{tiling_ifs}
Michael~F. Barnsley, Louisa~F. Barnsley, and Andrew Vince.
\newblock Tiling iterated function systems, 2020.

\bibitem{coding}
Douglas Lind and Brian Marcus.
\newblock {\em An Introduction to Symbolic Dynamics and Coding}.
\newblock Cambridge Mathematical Library. Cambridge University Press, second edition, 2021.

\bibitem{gifs}
R.~Mauldin and Stanley Williams.
\newblock Hausdorff dimension in graph directed constructions.
\newblock {\em Transactions of The American Mathematical Society}, 309:811--829, Feb 1988.

\end{thebibliography}
\end{document}